\documentclass[11pt]{amsart}

\usepackage{amsmath,amssymb,amscd,amsfonts,verbatim}

\usepackage[mathscr]{eucal}

\newcommand{\pionealg}{\pi_1^{\mathrm{alg}}}
\newtheorem{thm}{Theorem}[section]
\newtheorem{lem}[thm]{Lemma}

\newtheorem{prop}[thm]{Proposition}
\newtheorem{cor}[thm]{Corollary}

\newtheorem{assu-nota}[thm]{Assumption--Notation}
\theoremstyle{remark}

\newcommand{\into}{\hookrightarrow}

\newcommand{\Z}{\mathbb Z}

\newcommand{\pp}{\mathbb P}
\DeclareMathOperator{\Aut}{Aut}
\DeclareMathOperator{\Pic}{Pic}

\DeclareMathOperator{\Tr}{Tr}

\DeclareMathOperator{\Proj}{Proj}

\newcommand{\Si}{\Sigma}
\newcommand{\si}{\sigma}

\newcommand{\OO}{\mathcal{O}}

\numberwithin{equation}{section}

\title [ On the algebraic fundamental group of...]
{On the algebraic fundamental group of surfaces with $K^2\le3\chi$
}
\author{ Margarida Mendes Lopes \and Rita Pardini}

\thanks{The first author is a member of the Center for Mathematical
Analysis, Geometry and Dynamical Systems, IST TULisbon,  and the second author is a member
of G.N.S.A.G.A.-I.N.d.A.M. This research was partially supported by the
Italian project ``Geometria sulle variet\`a algebriche" (PRIN COFIN 2004)
and by FCT (Portugal) through program POCTI/FEDER and Project
POCTI/MAT/44068/2002.}

\begin{document}
\begin{abstract} Let $S$ be a minimal complex surface of general type with $q(S)=0$. We prove the following statements concerning the algebraic fundamental
group $\pionealg(S)$:
\begin{itemize}
\item Assume that $K^2_S\le3\chi(S)$. Then $S$ has an irregular \'etale
cover if and only if $S$ has a free pencil of hyperelliptic curves of
genus~$3$ with at least $4$ double fibres. 
\item If $K^2_S=3$ and $\chi(S)=1$, then $S$ has no irregular \'etale cover.
\item If $K^2_S<3\chi(S)$ and $S$ does not have any irregular \'etale cover,
then $|\pionealg(S)|\le9$. If $|\pionealg(S)|= 9$, then $K^2_S=2$, $\chi(S)=1$.
\end{itemize}
\medskip

\end{abstract}
\maketitle

\section { Introduction}

Every minimal surface $S$ of general type satisfies the Noether inequality:
$$K^2_S\ge 2\chi(S)-6.$$

It has been clear for a long time that the closer a surface is to the
Noether line $K^2=2\chi-6$, the simpler its algebraic fundamental group
is. In fact, Reid has conjectured that for $K^2<4\chi$ the algebraic fundamental group of $S$ is either finite or  it  coincides, up to finite group extensions, with the fundamental group of a curve of genus $g\geq 1$, i.~e. it is {\sl commensurable} with the fundamental group of a curve, (\cite[Conjecture 4] {miles1}, see also \cite{bpv}, pp. 294). 

In the case of irregular surfaces  or of regular surfaces having an irre\-gular \'etale cover,  Reid's conjecture follows from the Severi inequality, recently proved in \cite{severi},  which states that the Albanese map of an irregular surface with $K^2<4\chi$ is a pencil.

Indeed, let $S$ be  an irregular surface satisfying $K^2<4\chi$, let $a\colon S\to B$ be the Albanese pencil of $S$ and $F$ a
general fibre of $a$. The inclusion $F\into S$ induces a map
$\psi\colon \pionealg(F)\to\pionealg(S)$.
By \cite[Theorem~1]{xiaoslope} the image $H$ of $\psi$ is either 0 or $\Z_2$,
and $H=\Z_2$ is possible only if $F$ is hyperelliptic. The cokernel of $\psi$ is the so-called  {\em orbifold fundamental group} of the fibration $a$ (cf. \cite{cko},  \cite[Lemma 4.2]{catanese}).
If $a$ has no multiple fibres, then we have  an exact sequence:
 \begin{equation}\label{pisequence}
1\to H\to\pionealg(S)\to \pionealg(B)\to 1.
\end{equation}
If $a$ has multiple fibres, then it is possible to find a Galois cover  $B'\to B$  such that the fibration $a'\colon S'\to B'$ obtained from $a$ by base change and normalization has no multiple fibres and the map $S'\to S$ is \'etale. 
Since $\pionealg(S')$ is a normal subgroup  of $\pionealg(S)$ of  finite index, it follows that in any case the algebraic fundamental group of an irregular surface  satisfying $K^2<4\chi$  is commensurable with the fundamental group of a curve. Of course  the same is true for a regular surface satisfying $K^2<4\chi$ and having an irregular \'etale cover.

\medskip

 Reid's conjecture is still  open for surfaces not having an irregular cover. 
 However for surfaces satisfying $K^2<3\chi$ not only Reid's conjecture is true (\cite{miles1} and \cite{horikawaV}) but   work by
 several authors gives more precise results on the algebraic fundamental group  (cf. \cite{bombieri}, \cite{horikawaV}, \cite{miles1}, \cite{milesk2}, \cite{xiaohyp}, \cite{xiaoslope}). The picture that emerges from their work is the following:
\begin{itemize}
\item[(I)] if $K^2_S<2\chi(S)$, then $S$ is regular and $\pionealg(S)$ is finite.
\item[(II)] if $K^2_S<\frac{8}{3}\chi(S)$ and $S$ is irregular, then the Albanese map of $S$ is a pencil of curves of genus~2. If $K^2_S<\frac{8}{3}\chi(S)$ and $S$ is regular, then $\pionealg(S)$ is finite.
\item[(III)] if $K^2_S<3\chi(S)$ and $S$ is irregular, then the Albanese map of $S$ is a pencil of hyperelliptic curves of genus~2 or~3. If $S$ is regular, then either $\pionealg(S)$ is finite or there exists an irregular \'etale cover $X\to S$. The Albanese map of $X$ is a pencil of hyperelliptic curves of genus~3, which induces on $S$ a free pencil of hyperelliptic curves of genus~3 with at least $4$ double fibres. 
Conversely, if $S$ has such a pencil, then it admits an irregular \'etale cover.
\end{itemize}

These results give a good understanding of the algebraic fundamental group of a surface $S$ with $K^2<3\chi$ and infinite $\pionealg(S)$. 

In fact, if $S$ is irregular  and the Albanese map $a\colon S\to B$ has multiple fibres, then by statement (III) and by the adjunction formula we have $g=3$ and the multiple fibres are double fibres. Then there is a Galois cover $B'\to B$ with Galois group $G$ such that 
the $G$-cover $S'\to S$ obtained by base change and normalization is
\'etale and the induced fibration $a'\colon S'\to B'$ has no multiple
fibres. One can show that $G$ can be chosen to be a quotient of the
dihedral group of order 8. 
So we have an exact sequence:
$$1\to\pionealg(S')\to\pionealg(S)\to G\to 1.$$
and the group  $\pionealg(S')$ is described by sequence (\ref{pisequence}).

If $S$ is a regular surface such that $K^2_S<3\chi(S)$ and $\pionealg(S)$
is infinite, then using (III), one constructs an irregular \'etale Galois
cover $X\to S$ with Galois group $\Z_2$ or $\Z_2^2$ whose Albanese map is a pencil of curves of genus~3 without multiple fibres (more precisely, we have $\Z_2$ if the number $k$ of double fibres of $a$ is even and $\Z_2^2$ if $k$ is odd). Then the group $\pionealg(X)$ is a normal subgroup of $\pionealg(S)$ of index $2$ or $4$ which can be described as explained above. 

\medskip

However if the algebraic fundamental group of $S$ is finite then the above results give no additional information. \medskip

In this paper we give two improvements of the above results. 

We first extend part of (III) to surfaces on the line $K^2= 3\chi$:
\begin{thm}\label{main} Let $S$ be a minimal complex surface of general type with $q(S)=0$ and $K^2_S\le3\chi(S)$.

Then $S$ has an irregular \'etale cover if and only if there exists a fibration $f\colon S\to\pp^1$ such that:
\begin{enumerate}
\item the general fibre $F$ of $f$ is hyperelliptic of genus~$3$;
\item $f$ has at least $4$ double fibres.
\end{enumerate}
\end{thm}
This improvement is made possible by the Severi inequality.

In the case $p_g(S)=0$, Theorem~\ref{main} can be made more precise:
\begin{thm}\label{K3}
Let $S$ be a smooth minimal surface of general type with $p_g(S)=0$, $K^2_S=3$.

Then $S$ has no irregular \'etale cover.
\end{thm}
Theorem~\ref{K3} is sharp in a sense, since there are examples, due to Keum
and Naie (cf. \cite{naie}), of surfaces with $K^2=4$ and $p_g=0$ that have
an irregular cover. 

On the other hand, it remains an open question whether
the algebraic fundamental group of a surface with $K^2=3$ and $p_g=0$ is finite or more generally whether the algebraic fundamental
group of a surface with $K^2=3\chi$ that has no \'etale irregular cover
is finite.  

In even greater generality  one would like to know whether the algebraic fundamental
group of a surface with $K^2<4\chi$ that has no \'etale irregular cover
is finite, deciding thus  Reid's conjecture. This is a very challenging problem, which however does not seem possible to resolve with the methods of the present paper. 

\medskip 

Finally, 
we bound the cardinality of $\pionealg(S)$ in the case when it is a finite group:
\begin{thm}\label{T2}
Let $S$ be a minimal surface of general type such that $K^2_S<3\chi(S)$.
If $S$ has no irregular \'etale cover, then $\pionealg(S)$ is a finite group of order $\le9$.

Moreover, if $\pionealg(S)$ has order 9, then $\chi(S)=1$ and $K^2_S=2$,
namely $S$ is a numerical Campedelli surface.
\end{thm}
This bound is sharp, since there are examples of surfaces with $p_g=0$, $K^2=2$ and $\pionealg(S)=\Z_9$, $\Z_3^2$ (cf. \cite[Ex. 4.11]{xiaog2}, \cite{pi9}).  

By this  theorem  only a very short list of finite groups can occur as the algebraic fundamental groups of surfaces with $K^2\leq 3\chi-1$.  The list is even more restricted if $K^2\leq 3\chi-2$: in \cite{3chi-2} it is shown that in this case $|\pionealg(S)|\le 5$, with equality holding only for surfaces with $K^2_S=1$ and $p_g(S)=0$. Moreover $|\pionealg(S)|=3$  is possible only for $2\le \chi(S)\le 4$ and $K^2=3\chi-3$.
\medskip

{\bf Notation and conventions.} We work over the complex numbers. All varieties are projective algebraic.
We denote by $\chi$ or $\chi(S)$ the holomorphic Euler characteristic of the structure sheaf of the surface $S$.

\section {The proof of Theorem~\ref{main}}

In this section we assume that $S$ is a minimal complex surface of general type with $q(S)=0$ and $K^2_S\le3\chi(S)$. In order to prove Theorem~\ref{main} we need some intermediate steps.

\begin{lem}\label{fibre3}
Let $\rho\colon Z\to S$ be an \'etale cover such that $q(Z)>0$. 

Then the Albanese pencil $a\colon Z\to A$ induces a fibration $f\colon S\to\pp^1$ such that:
\begin{enumerate}
\item the general fibre $F$ of $f$ is a curve of genus~$3$;

\item $f$ has at least 4 double fibres.
\end{enumerate}
Moreover, all irregular \'etale covers of $S$ induce the same fibration
$f\colon S\to\pp^1$.
\end{lem}
\begin{proof}

If $\rho \colon Z\to S$ is an irregular \'etale cover, then the Galois closure of $\rho$ is an irregular Galois \'etale cover. We denote by $\pi\colon Y\to S$ 
a minimal element of the set of irregular Galois \'etale covers of $S$.

Denote by $d$ the degree of $\pi$. The surface $S$ is minimal of general type with $K^2_Y=dK^2_S$, $\chi(Y)=d\chi(S)$.
Hence $K^2_Y\le3\chi(Y)<4\chi(Y)$ and therefore, by the Severi inequality
(\cite{severi}), the image of the Albanese map of $Y$ is a curve. Write
$a\colon Y\to B$ for the Albanese pencil, and let $b$ be the genus of $B$
and $g$ the genus of the general fibre $F$ of $a$. The Galois group $G$
of $\pi$ acts on the curve $B$. This action is effective by the assumption
that $\pi$ is minimal among the irregular \'etale covers of $Y$.
Hence we have a commutative diagram:
\begin{equation} \label{diagram}
\begin{CD}
Y @>{\pi}>>S \\
@V{a}VV @VV{f}V\\
B @>{\bar{\pi}}>> \pp^1 
\end{CD}
\end{equation}

The map $\bar{\pi}$ is a Galois cover with group $G$ and the general fibre of $f$ is also equal to $F$. Since the map $\pi$ is obtained from $f$ by taking base change with $\bar{\pi}$ and normalizing, the fibre of $f$ over a point $x$ of $\pp^1$ has multiplicity equal to the ramification order of $\bar{\pi}$ over $x$. Notice that, since $\pp^1$ is simply connected, the branch divisor of $\bar{\pi}$ is nonempty and therefore the fibration $f$ always has multiple fibres. Notice also that, since $S$ is of general type, the existence of multiple fibres implies $g\ge 3$.

We remark that the fibration $a$ is not smooth and isotrivial. In fact, if this were the case then $Y$ would be a free quotient of a product of curves, hence it would satisfy $K^2_Y=8\chi(Y)$.
Hence we may define the slope of $a$ (cf. \cite{xiaoslope}): $$\lambda(a):=\frac{K^2_Y-8(b-1)(g-1)}{\chi(Y)-(b-1)(g-1)}.$$ The slope inequality (\cite{xiaoslope}, cf. also \cite{ch}, \cite{stoppino}) gives 
\begin{equation}\label{slope}
4(g-1)/g\le \lambda(a)\le K^2_Y/\chi(Y)=K^2_S/\chi(S)\le3 \end{equation}
where the second inequality is a consequence of $b>0$.
Hence we get $g=3$ or $g=4$.

Assume $g=4$. In this case (\ref{slope}) becomes:
$$3\le \lambda(a)\le K^2_S/\chi(S)\le3.$$ It follows that the slope inequality is sharp in this case and $K^2_S=3\chi(S)$. By \cite[Prop. 2.6]{konnoslope}, this implies that $F$ is hyperelliptic.
Let $\si$ be the involution of $S$ induced by the hyperelliptic involution on the fibres of $f$. The divisorial part $R$ of the fixed locus of $\si$ satisfies $FR=10$.
As remarked above, $f$ has at least a fibre of multiplicity $m>1$, that
we denote by $mA$. 
Since $g=4$, by the adjunction formula $\frac{6}{m}$ is divisible by $2$,
yielding $m=3$. Hence $3AR=10$, a contradiction. So we have proved $g=3$.

Using the adjunction formula again, we see that the multiple fibres of $f$ are double fibres, hence all the branch points of $\bar{\pi}$ have ramification order equal to $2$. Let $k$ be the number of branch points of $\bar{\pi}$. By applying the Hurwitz formula to $\bar{\pi}$, we get $k\ge 4$.

Given an irregular \'etale cover $\rho\colon Z\to S$, we can always find an \'etale cover $W\to S$ which dominates both $Z$ and $Y$. The Albanese pencil of $W$ is a pullback both from $Y$ and from $Z$, hence the fibrations induced on $S$ by the Albanese pencils of $Z$, $W$ and $Y$ are the same. 
\end{proof}

We introduce some more notation. Assume that $f\colon S\to \pp^1$ is the
fibration defined in Lemma \ref{fibre3}. Let $\bar{\pi}\colon B\to \pp^1$
be the double cover branched on 4 points corresponding to double fibres
$2F_1,\dots,2F_4$ of $f$ and $\pi\colon Y\to S$ the \'etale double cover
obtained by base change with $\bar{\pi}$ and normalization, as in diagram
(\ref{diagram}). Then $K^2_Y=2K^2_S$, $\chi(Y)=2\chi(S)$ and $q(Y)=1$. We
write $\eta:=F_1+F_2-F_3-F_4$. Clearly, $\eta$ has order~2 in $\Pic(S)$
and $\pi$ is the \'etale double cover corresponding to $\eta$.

\begin{lem} \label{hyp} The general fibre $F$ of $f$ is hyperelliptic.
\end{lem}
\begin{proof}

Assume by contradiction that $F$ is not hyperelliptic and consider the pencil $a\colon Y\to B$, whose general fibre is also equal to $F$. Set $\mathcal E:=a_*\omega_Y$ and denote by $\psi\colon Y\to \pp(\mathcal E)$ the relative canonical map, which is a morphism by Remark 2.4 of \cite{konnoslope}. 
Let $V$ be the image of $\psi$.
The surface $V$ is a relative quartic in $\pp(\mathcal E)$ and, by Lemma 3.1 and Theorem~3.2 of \cite{konnoslope}, its singularities are at most rational double points. The map $\psi$ is birational and it contracts precisely the nodal curves of $Y$, which are all vertical since $B$ has genus~1. Hence $V$ is the
canonical model of $Y$.

Let $\iota$ be the involution associated to the cover $Y\to S$. This involution induces automorphisms of $B$, $\mathcal E$, $\pp(E)$ and $V$ (that we denote again by $\iota$) compatible with $a$, $\psi$ and the inclusion $V\subset \pp(\mathcal E)$. Given $b\in B$, 
write $\pp^2_b$ for the fiber of $\pp(\mathcal E)$ over $b$ and $V_b:=V\cap \pp^2_b$. The curve $V_b$ is a plane quartic inside $\pp^2_b$. For every $b\in B$, the map $\iota$ induces a projective isomorphism between $\pp^2_b$ and $\pp^2_{\iota(b)}$ that restricts to an isomorphism of $V_b$ with $V_{\iota(b)}$. In particular, if $b$ is one of the four fixed points of $\iota$ on $B$, then $\iota$ induces an involution of $\pp^2_b$ that preserves the quartic $V_b$. Since the fixed locus of an involution of the plane contains a line, it follows that $\iota$ has at least a fixed point on $V_b$. In particular, the action of $\iota$ on $V$ is not free.

On the other hand, one checks that a fixed point free auto\-morphism of
a minimal surface of general type induces a fixed point free automorphism
of the canonical model. So we have a contradiction.
\end{proof}

We can now give:
\begin{proof}[Proof of Theorem~\ref{main}]
The ``if'' part is a consequence of Lemma \ref{fibre3} and Lemma \ref{hyp}.
Conversely, if $S$ has a fibration with 4 double fibres $2F_1,\dots,2F_4$
then the \'etale double cover associated with $\eta:=F_1+F_2-F_3-F_4$ has irregularity equal to 1.
\end{proof}

\section {The proof of Theorem~\ref{K3}}

In this section we let $S$ denote a
smooth minimal surface of general type with $p_g(S)=0$, $K^2_S=3$.
To prove Theorem~\ref{K3} we argue by contradiction.

Thus assume that $S$ has an irregular \'etale cover. Then by
Theorem~\ref{main} there exists a fibration $f\colon S\to\pp^1$ whose
general fibre is hyperelliptic of genus~3 and with at least 4 double
fibres $2F_1,\dots,2F_4$. As before,
denote by $\pi\colon Y\to S$ the \'etale double cover given by $\eta=F_1+F_2-F_3-F_4$ and by $\iota$ the involution associated with $\pi$. The invariants of $Y$ are: $q(Y)=1$, $p_g(Y)=2$, $K^2_Y=6$.

The hyperelliptic involution on the fibres of $a\colon Y\to B$ and $f\colon S\to Y$ induces involutions $\tau$ of $Y$ and $\si$ of $S$. By construction, these involutions are compatible with the map $\pi\colon Y\to S$, namely we have $\pi\circ \tau=\si\circ\pi$.
We denote by $p\colon S\to \Si:=S/\si$ the quotient map.

\begin{lem}\label{commute}
The involutions $\tau$ and $\iota$ of $Y$ commute.
\end{lem}
\begin{proof}
Denote by $h$ the composite map $Y\to S\to\Si$. By construction, both $\iota$ and $\tau$ belong to the Galois group $G$ of $h$. Since $h$ has degree 4 and $\iota$ and $\tau$ are involutions, the group $G$ is isomorphic to $\Z_2\times \Z_2$ and $\iota$ and $\tau$ commute.
\end{proof}

\begin{lem}\label{16}
The involution $\iota\tau$ has at least $16$ isolated fixed points on $Y$.
\end{lem}
\begin{proof}
Let $q\colon Y\to Z:=Y/\iota\tau$ be the quotient map. The surface $Z$ is
nodal. The regular $1$-forms and $2$-forms of $Z$ correspond to the
elements of $H^0(Y, \Omega^1_Y)$, respectively $H^0(Y,\omega_Y)$, that
are invariant under the action of $\iota\tau$. By the same argument,
since $p_g(S)=p_g(Y/\tau)=0$, both $\iota$ and $\tau$ act on
$H^0(Y,\omega_Y)$ as multiplication by $-1$. It follows that $\iota\tau$
acts trivially on $H^0(Y,\omega_Y)$ and $p_g(Z)=2$. Since $\iota$ acts
on $B$ as an involution with quotient $\pp^1$ and $\tau$ acts trivially
on $B$, it follows that the action of $\iota\tau$ on $B$ is equal to
the action of $\iota$ and that $q(Z)=0$.

Let $D$ be the divisorial part of the fixed locus of $\iota\tau$ on $Y$ and let $k$ be the number of isolated fixed points of $\iota\tau$. We recall the Holomorphic Fixed Point formula (see \cite{as}, p.566):
$$\sum_i(-1)^i \Tr(\iota\tau|H^i(Y,\OO_Y))= (k-K_YD)/4.$$
By the above considerations, this can be rewritten as:
$$k=16+K_YD.$$
The statement now follows from the fact that $K_Y$ is nef.
\end{proof}

\bigskip

\begin{proof}[Proof of Theorem~\ref{K3}]
By Lemma \ref{commute}, the involution $\iota\tau$ of $Y$ induces $\si$ on $S$. By Lemma \ref{16}, $\iota\tau$ has at least 16 isolated fixed points. Since the images on $S$ of these points are isolated fixed points of $\si$, the involution $\si$ has at least 8 isolated fixed points.
On the other hand, by \cite[Prop. 3.3]{CCM} there are at most $K^2_S+4=7$ isolated fixed points of $\si$.
So we have a contradiction, and thus $S$ has no irregular \'etale cover.
\end{proof}

\section{The proof of Theorem~\ref{T2}}

To prove Theorem~\ref{T2} we will use the following two results proved
in \cite[Cor. 5.8]{beauville}, although not stated explicitly.

\begin{prop} \label{be} Let $Y$ be a surface of general type such that
the canonical map of $Y$ has degree~$2$ onto a rational surface. If $G$
is a group that acts freely on $Y$, then $G=\Z_2^r$, for some $r$.
\end{prop}
\begin{proof} The group $G$ is finite, since a surface of general type has finitely many automorphisms.

Let $T$ be the quotient of $Y$ by the canonical involution. The surface $T$ is rational, with canonical singularities, and $G$ acts on $T$.

 Since $T$ is rational, each element $g\in G$ acts with fixed points. The argument in the proof of \cite[Cor. 5.8]{beauville} shows that each $g$ has order $2$, hence $G=\Z_2^r$.
\end{proof}

\begin{cor}\label{Z2}
Let $S$ be a minimal surface of general type such that $K^2_S<3\chi(S)$,
and  $S$ has no irregular \'etale cover. If   $Y\to S$ is an \'etale G-cover, then either  $|G|\leq 10$ or $G=\Z_2^r$, for some $r\geq 4$.
\end{cor}
\begin{proof}
Let $\pi\colon Y\to S$ be an \'etale G-cover of degree $d>10$.
By assumption we have $q(Y)=0$ and 
$K^2_Y<3p_g(Y)-7$, and therefore
the canonical map of $Y$ is 2-to-1 onto a rational surface by
\cite[Theorem~5.5]{beauville}. Hence $G=\Z_2^r$ for some $r\ge 4$ by Proposition \ref{be}. 
\end{proof} 

For related statements see the results of \cite{xiaohyp} on hyperelliptic surfaces and the results of \cite{ak} and \cite{konnopg}.

 We remark that the next result is well known for the cases $\chi(S)=1$ and $K^2_S=1$ or $2$ (\cite{milesk2}).
\begin{prop}\label{pi9}
Let $S$ be a minimal surface of general type with $K^2_S<3\chi(S)$.
If $S$ has no irregular \'etale cover, then $|\pionealg(S)|\le9$.\par

\end{prop}
\begin{proof}

Let $Y\to S$ be an \'etale G-cover.
By Corollary \ref{Z2}, it is enough to exclude the following possibilities: a) $G=\Z_2^r$ for some $r\ge 4$, and b) $|G|=10$.

Consider case a) and assume by contradiction that $\pi\colon Y\to S$ is a Galois \'etale cover with Galois group $G=\Z_2^4$. By \cite{miyaoka1}, $\chi(S)\geq 2$. 
We have $\chi(Y)=16\chi(S)\ge 32$ and $K^2_Y<3(\chi(Y)-5)$. Notice that,
since $K^2_Y<3\chi(Y)-10$, by \cite[Theorem~5.5]{beauville} the surface
$Y$ has a pencil of hyperelliptic curves. Hence $Y$ satisfies the
assumptions of \cite[Theorem~1]{xiaohyp} and there exists a unique free
pencil $|F|$ of hyperelliptic curves of genus $g\le3$ on $Y$. The action
of $G$ preserves $|F|$ by the uniqueness of $|F|$. Since $\Aut(\pp^1)$
does not contain a subgroup isomorphic to $\Z_2^3$, there is a subgroup
$H<G$ of order $\ge 4$ that maps every curve of $|F|$ to itself.
Since the action of $G$ on $Y$ is free, this implies that $g-1$ is divisible by 4, contradicting $g\le3$ and $S$ of general type.

Consider now case b) and assume by contradiction that $\pi \colon Y\to S$ is a Galois cover with Galois group $G$ of order 10. For $K^2_S<3\chi(S)-1$, we have $K^2_Y<3\chi(Y)-10$ and, as in the proof of Corollary \ref{Z2}, $G$ is of the form $\Z_2^a$, a contradiction. So we have $K^2_S=3\chi(S)-1$, $K^2_Y=3\chi(Y)-10$, $q(Y)=0$ and so, by \cite{ak}, the canonical map of $Y$ is either birational or 2-to-1 onto a rational surface. By Proposition \ref{be}, the last possibility does not occur, since $G$ has order 10.

The surface $Y$ satisfies $p_g(Y)=10\chi(S)-1\ge 9$. Surfaces on the Castelnuovo line $K^2=3\chi-10$ with birational canonical map are classified (cf. \cite{harris}, \cite{rick} and \cite{ak}): for $p_g(Y)\ge 8$, the canonical model
$V$ of $Y$ is a relative quartic inside a $\pp^2$-bundle
\[
\pp:=\Proj(\OO_{\pp^1}(a)\oplus \OO_{\pp^1}(b)\oplus\OO_{\pp^1}(c)),
\]
where $0\le a\le b\le c$ and $a+b+c=p_g(Y)+3$. 

If the Galois group $G$ preserves the fibration $f\colon V\to\pp^1$
induced by the projection $\pp\to\pp^1$, then, as in Lemma~\ref{hyp},
we obtain a contradiction by considering the action on $V$ of an
element of order~2 of $G$. 

So, to conclude the proof we just have to show that $G$ preserves $f$. Let $W$ be the image of $\pp$ via the tautological linear system. By the results of \cite{ak}, \cite{harris}, \cite{rick}, the threefold $W$ is the intersection of all the quadrics that contain the canonical image of $Y$ and therefore it is preserved by the automorphisms of $V$. One checks that $W$ has a unique ruling by
planes which induces the fibration $f$ on $V$. Therefore every automorphism
of $V$ preserves the fibration $f$. \end{proof}

To obtain the statement of Theorem~\ref{T2} we now show the following:

\begin{prop}\label{campedelli}
Let $S$ be a minimal surface of general type with $K^2_S<3\chi(S)$. If $|\pionealg(S)|=9$, then $\chi(S)=1$ and $K^2_S=2$, namely $S$ is a numerical Campedelli surface.
\end{prop}

\begin{proof}
Suppose that $|\pionealg(S)|=9$ and $\chi(S)\geq 2$. The argument in the
proof of Proposition~\ref{pi9} shows that $K^2_S=3\chi(S)-1$. Let
$\pi\colon Y\to S$ be the universal cover. We have $K^2_Y=3p_g(Y)-6$, $p_g(Y)=9\chi(Y)-1\ge 17$. By \cite[Lem. 2.2]{konnopg} the bicanonical map of $Y$ has degree 1 or 2. Arguing as in the proof of Proposition~\ref{pi9}, one shows that the bicanonical map of $Y$ is birational. Then, since $p_g(Y)\ge 11$, by the results of \cite{konnopg} the situation is analogous to the case of a surface with $K^2=3p_g-7$ and birational canonical map. Namely, the intersection of all the quadrics through the canonical image of $Y$ is a threefold $W$, which is the image of a $\pp^2$-bundle $\pp:=\Proj(\OO_{\pp^1}(a)\oplus\OO_{\pp^1}(b)\oplus\OO_{\pp^1}(c))$ via the tautological linear system, and $Y$ is birational to a relative quartic of $\pp$.
In particular, there is a fibration $f\colon Y\to\pp^1$ with general fibre
a nonhyperelliptic curve of genus~3. One can show as above that the Galois
group $G=\pionealg(S)$ of $\pi$ preserves $f$. Then we obtain a
contradiction, since the multiple fibres of a genus~3 fibration are double fibres and a smooth genus~3 curve does not admit a free action of a group of order 9.
\end{proof}

{\em Remark.} Numerical Campedelli surfaces with fundamental group $\Z_9$ and $\Z_3^2$ do exist (cf. \cite[Ex. 4.11]{xiaog2}, \cite{pi9}).

\bigskip

\bigskip

\begin{minipage}{13cm}
\parbox[t]{6.7cm}{Margarida Mendes Lopes\\
Departamento de Matem\'atica\\
Instituto Superior T\'ecnico\\
Universidade T{\'e}cnica de Lisboa\\
Av.~Rovisco Pais\\
1049-001 Lisboa, PORTUGAL\\
mmlopes@math.ist.utl.pt
} \hfill
\parbox[t]{5.5cm}{Rita Pardini\\
Dipartimento di Matematica\\
Universit\`a di Pisa\\
Largo B. Pontecorvo, 5\\
56127 Pisa, Italy\\
pardini@dm.unipi.it}
\end{minipage}

\end{document}